\documentclass[11pt,a4paper,subeqn,reqno]{article}

\usepackage{amsmath,amsthm}
\usepackage{amssymb}
\usepackage[mathscr]{eucal}
\usepackage[all]{xy}
\usepackage{mathrsfs}
\usepackage{hyperref}
\usepackage{appendix}
\usepackage{authblk}
\usepackage{graphicx}
\usepackage[misc]{ifsym}
\usepackage{tikz,tikz-cd}
\usepackage{url}
\usepackage{hyperref}

\allowdisplaybreaks 

\usepackage[text={165mm,250mm},centering]{geometry}

\newtheorem{theorem}{Theorem}[section]
\newtheorem{lemma}[theorem]{Lemma}

\newtheorem{corollary}[theorem]{Corollary}
\newtheorem{proposition}[theorem]{Proposition}

\theoremstyle{definition}
\newtheorem{example}[theorem]{Example}
\newtheorem{definition}[theorem]{Definition}
\newtheorem{definition-lemma}[theorem]{Definition-Lemma}
\newtheorem{definition-theorem}[theorem]{Definition-Theorem}

\newtheorem{remark}[theorem]{Remark}

\newtheorem{proof*}{Proof}

\begin{document}

\title{Langton's type Theorem on Algebraic Orbifolds}

\author[1]{Yonghong Huang\thanks{Email: huangyh329@mail.sysu.edu.cn (yonghonghuangmath@gmail.com)}}


\renewcommand\Affilfont{\small}

\affil[1]{School of Mathematics, Sun Yat-sen University, Guangzhou 510275, P. R. China}

\date{}

\maketitle
\renewcommand{\abstractname}{Abstract}

\begin{abstract}
In this paper, we prove the Langton's type theorem on separatedness and properness for the moduli functor of torsion free semistable sheaves on algebraic orbifolds over an algebraically closed field $k$.
\end{abstract}

\section{Introduction}
\hspace*{1em}Let $\mathcal X$ be a smooth algebraic orbifold (Def.\ref{def al orb} and Remark \ref{remark alg orb}) over an algebraically closed field $k$. We consider the moduli functor $\mathcal M$ of modified slope (Def.\ref{def mslope}) semistable torsion free sheaves on $\mathcal X$. Following \cite{hl}, \cite{lan}, \cite{sp}, we define the following  functor:
\[
\widehat{\mathcal M}: ({\rm Sch}/k)^o\longrightarrow (\rm Sets)
\]
as follows. Let $T$ be a $k$-scheme  and let $\widehat{\mathcal M}(T)$ be the set of isomorphism classes of $T$-flat families of torsion free semistable sheaves on $\mathcal X$. If $f:T^\prime\rightarrow T$ is a morphism of schemes, let $\widehat{\mathcal M}(f)$ be the morphism obtained by pulling back sheaves via the morphism $f_{\mathcal X}=\text{id}_{\mathcal X}\times f$, i.e \[
\widehat{\mathcal M}(T)\longrightarrow\widehat{\mathcal M}(T^\prime),\quad [E]\longmapsto [f^*_{\mathcal X}E].
\]
Then, the moduli functor $\mathcal M$ is defined to be the quotient functor of $\widehat{\mathcal M}$ by equivalence relation $\sim$:
\[
E\sim E^\prime,\quad\text{for $E,E^\prime\in\mathcal M^\prime(T)$ if and only if there is a line bundle $L$ on $T$ such that $E^\prime=p_2^*L\otimes E$},
\]
where $p_2:\mathcal X\times T\rightarrow T$ is the projection onto $T$. In general, the moduli functor $\mathcal M$ is not representable. In fact, if $\mathcal X$ is a projective scheme and there is a properly semistable sheaf on $\mathcal X$, then the moduli functor $\mathcal M$ can not be represented (Lemma 4.1.2 in \cite{hl}). In the case that $\mathcal M$ is representable, Nironi has shown that the corresponding moduli scheme is proper over $k$ (Theorem 6.22 in \cite{fn}). But, by the Grothendieck's valuative criteria, we can also consider the separatedness and properness of $\mathcal M$ directly. Indeed, Langton \cite{lan} has showed that the moduli functor of slope semistable torsion free sheaves on smooth projective varieties over $k$ is separated and proper. Maruyama \cite{mm}, Mehta and Ramanathan \cite{mr} generalised Langton'results to Gieseker stability. In recent years, many problems about the modui functor of semistable sheaves on algebraic orbifolds are concerned. There is not a similar result on algebraic orbifolds. For the researchers' convenience, we generalize the result of Langton and prove that the moduli functor of slope semistable torsion free sheaves on algebraic orbifolds is separated and proper. For the case of Gieseker stability, the similar result can be obtain following the line of Maruyama\cite{mm}, Mehta and Ramanathan \cite{mr}, for the sake of the key Lemma \ref{lemm ext 1} on algebraic orbifolds (which corresponds to the Proposition 6 in \cite{lan}). In the next paragraph, we give the precise description of the problem. Let $R$ be a discrete valuation ring over $k$ with maximal ideal $(\pi)$ and residue field $k$. The quotient field of $R$ is $K$. Consider the following cartesian diagram:
\begin{equation*}
\xymatrix{
\mathcal X_K \ar@{^(->}[r]^i \ar[d] & \mathcal X_R \ar[d] & \mathcal X_k \ar[d] \ar@{_(->}[l]_{j}\\
     \text{Spec}(K) \ar@{^(->}[r] & \text{Spec}(R) & \text{Spec}(k) \ar@{_(->}[l] }
\end{equation*}
where $\mathcal X_R=\mathcal X\times \text{Spec}(R)$, $\mathcal X_K=\mathcal X\times \text{Spec}(K)$ and
$\mathcal X_k=\mathcal X\times\text{Spec}(k)=\mathcal X$.

Consequently, we have:
\begin{enumerate}
\item [(1)] $\mathcal M$ is separated if and only if two families $E_R$, $E_R^\prime$ of torsion free semistable sheaves over $\text{Spec}(R)$ agreeing on the generic fiber $\mathcal X_K$, then they agree on $\mathcal X_R$;
\item [(2)]$\mathcal M$ is proper if and only if every torsion free semistable sheaves $E_K$ on $\mathcal X_K$ can be uniquely extend to a flat family of torsion free sheaves on $\mathcal X_R$, under isomorphism.
\end{enumerate}
We state our main results:
\begin{theorem}
Assume that $E_K$ is a torsion free sheaf on $\mathcal X_K$. Then
\begin{enumerate}
\item If $E_1$ and $E_2$ are two coherent subsheaves of $i_*E_K$ on $\mathcal X_R$ such that $i^*E_1=i^*E_2=E_K$ and $j^*E_1$, $j^*E_2$ are semistable torsion free sheaves on $\mathcal X_k$, at least one of which is stable, then there is an integer $p$ such that $E_1=\pi^p E_2$.
\item If $E_K$ is semistable, then there exists a coherent subsheaf $E\subseteqq {i_*}E_K$ such that $i^*E=E_K$ and $j^*E$ is torsion free and semistable on $\mathcal X_k$.
\end{enumerate}
\end{theorem}
As an application of the Theorem 1.1, in a forthcoming paper \cite{HJ}, we use it to show that the Hitchin map on the moduli space of Higgs bundles on Deligne-Mumford curves is proper.

\section{Torsion free sheaves on algebraic orbifolds}\label{sect tor sh}
Throughout this paper, we work over a fixed algebraically closed field $k$. All schemes, algebraic spaces and stacks and morphisms among them are of finite type. In the following, we recall some basic knowledge about torsion free sheaves on algebraic orbifolds. For more details, we refer the reader to \cite{aov}, \cite{dm}, \cite{ak}, \cite{fn} and \cite{vistoli}.

\begin{definition}[Tame Deligne-Mumford stacks]
Let $\mathcal X$ be a Deligne-Mumford stack with coarse moduli space $p:\mathcal X\rightarrow X$. Then, $\mathcal X$ is tame if the pushforward functor $p_*:\text{QCoh}(\mathcal X)\rightarrow \text{QCoh}(X)$ is exact, where $\text{QCoh}(-)$ is the category of quasicoherent sheaves.
\end{definition}

\begin{definition}[Algebraic Orbifolds]\label{def al orb}
Let $\mathcal{X}$ be a Deligne-Mumford tame stack over $k$, which is isomorphic to a separated global quotient $\big[Z/G\big]$, where $Z$ is an algebraic space over $k$ and $G$ is a subgroup scheme (a locally closed subscheme which is a subgroup) of some $\text{GL}_{N,k}$. If the generic stabilizer of $\mathcal X$ is trivial, then $\mathcal X$ is called an algebraic orbifold over $k$.
\end{definition}

\begin{definition}\label{def irr and int}
A Deligne-Mumford stack $\mathcal X$ is called irreducible if it is not the union of two proper closed subsets, where the closed sets in $\mathcal X$ means reduced closed substacks of $\mathcal X$. It is called integral if it is both irreducible and reduced.
\end{definition}
\begin{remark}\label{remk irreducible}
A Deligne-Mumford stack $\mathcal X$ is irreducible if and only if its coarse moduli space is irreducible. In fact, there is a bijection between the closed subsets of $\mathcal X$ and the closed subsets of $X$, as pointed out by Conrad in \cite{bc}.
\end{remark}

Nironi \cite{fn} introduces the notion of projective (quasi-projective) Delgine-Mumford stack:

\begin{definition}[Projective (quasi-projective) Deligne-Mumford stack]
Let $\mathcal X$ be a Deligne-Mumford stack over a field $k$. We say $\mathcal X$ is projective (quasi-projective) over $k$ if it is a tame separated global quotient with projective (quasi-projective) coarse moduli scheme.
\end{definition}

\begin{remark}\label{remark alg orb}
In this paper, we only consider $\mathcal{X}$ to be a projective algebraic orbifold, which is irreducible and its coarse moduli space is a projective scheme over $k$.
\end{remark}

\begin{example}[Weighted Projective line]
The weighted Projective lines $\mathbb{P}(n,m)$ are algebraic orbifolds, when $m$ and $n$ are coprime.
\end{example}
For more example, the reader can consult \cite{ak}.
As point out by Nironi, for a stack, there is no very ample invertible sheaves unless it is an algebraic space. However, under certain hypothesis, there exist locally free sheaves, called generating sheaves, which behave like very ample sheaves.

\begin{definition}[Generating sheaf]
Let $\mathcal{X}$ be a tame Deligne-Mumford stack and let $\pi:\mathcal{X}\rightarrow X$ be the coarse moduli space of $\mathcal{X}$. A locally free sheaf $\mathcal E$ of $\mathcal X$ is said to be a generating sheaf if for any
quasi-coherent sheaf $F$, the following map
\[
\pi^{\ast}(\pi_{\ast}({\mathcal E}^\vee\otimes F))\otimes\mathcal E\longrightarrow F
\]
is surjective.
\end{definition}
Olsson and Starr proved the existence of the generating sheaves. Also, the generating sheaf is stable for arbitrary base change on the coarse moduli space.
\begin{proposition}\cite{os}
\begin{enumerate}
\item Let $\mathcal X$ be a separated Deligne-Mumford tame stack which is a global quotient over $k$,  then there is a locally free sheaf $\mathcal E$ over $\mathcal X$ which is a generating sheaf for $\mathcal X$.
\item Let $\pi:\mathcal X\rightarrow X$ be the moduli space of $\mathcal X$ and $f:X^\prime\rightarrow X$ a morphism of algebraic spaces over $k$. Moreover, we have the following cartesian diagram:
    \[
    \xymatrix{
      \mathcal X^\prime \ar[d]_{p} \ar[r] & \mathcal X \ar[d]^{\pi} \\
       X^\prime \ar[r]^{f} & X  }
    \]
   and $p^*\mathcal E$ is a generating sheaf for $\mathcal X^\prime$.
\end{enumerate}
\end{proposition}
For a Deligne-Mumford stack $\mathcal X$ with projective coarse moduli scheme over a field of characteristic zero, the existence of the generating sheaf is equivalent to $\mathcal X$ is a global quotient stack.

\begin{proposition}\cite{ak}
For a Deligne-Mumford stack $\mathcal X$ over $k$ and $char k=0$, the following are equivalent.
 \begin{enumerate}
 \item $\mathcal X$ has a projective coarse moduli space and is a quotient stack.
 \item $\mathcal X$ has a projective coarse moduli space and possesses a generating sheaf.
 \item $\mathcal X$ can be embedded into a smooth Deligne-Mumford stack with projective
 coarse moduli space.
 \end{enumerate}
 \end{proposition}

As the case of schemes, the support of coherent sheaves on Deligne-Mumford stacks can be defined in the following way.

\begin{definition}[Support of Coherent sheaf]\label{def sp ch}
Let $\mathcal X$ be a Deligne-Mumford stack over $k$ and let $F$ be a coherent sheaf on $\mathcal X$. The support
${\rm supp}(F)$ of $F$ is the closed substack defined by the sheaf of ideals
\[
\xymatrix@C=0.5cm{
0 \ar[r] & \mathcal I_{F} \ar[r] & \mathcal{O}_{\mathcal X} \ar[r] &  {\mathscr{H}om}_{\mathcal{O}_{\mathcal X}}(F,F)}.
\]
\end{definition}

\begin{definition}[Torsion free sheaf]\label{def tfs}
Let $\mathcal X$ be an projective Deligne-Mumford stack over $k$. A coherent sheaf $F$ is said to
be a torsion free sheaf if for every nonzero subsheaf $G\subseteq F$, the dimension of $\text{supp}(G)$
is $\text{dim}\mathcal X$.
\end{definition}
The torsion freeness of a coherent sheaf on a Deligne-Mumford stack is equivalent to its restriction to an \'etale covering.(Remark 3.3 in \cite{fn}).
\begin{lemma}[\cite{fn}]
With the same hypothesis as above, $F$ is a torsion free sheaf if and only if there is an \'etale covering
$f:U\rightarrow\mathcal X$ such that the restriction of $F$ to $U$ is torsion free.
\end{lemma}

\begin{proposition}\label{pro rk}
Assume that $\mathcal X$ is an integral projective Deligne-Mumford stack over $k$ and $F$ is a coherent sheaf on
$\mathcal X$. Then, there exists an open substack $\mathcal X^o$, such that the restriction $F|_{\mathcal X^o}$ to $\mathcal X^o$ of $F$ is locally free.
\end{proposition}
\begin{proof}
Take an \'etale covering $f:U\rightarrow\mathcal X$ such that $U$ is finite type over $k$. We have the following cartesian diagram:
\[
  \xymatrix{
 U\times_{\mathcal X}U \ar[d]_{\quad pr_1} \ar[r]^{\quad pr_2}
                & U \ar[d]^{f}  \\
  U \ar[r]_{f} & \mathcal X }
\]
Denote $U{\times}_{\mathcal X}U$ by $R$. Then $R\underset{t}{\overset{s}{\rightrightarrows}}U$
is an algebraic groupoid, where $s=pr_1$ and $t=pr_2$.
Denote $f^*F$ by $F^\prime$. By the $2$-commutativity of above diagram, there is an isomorphism
\[
\phi: s^{*}F^\prime\longrightarrow t^{*}F^\prime.
\]
Because $U$ is reduced, there exists unique maximal nonempty open subset $U^\prime\subset U$
such that $F^\prime|_{U^\prime}$ is locally free. By the flatness of morphism $s$, $s^{-1}(U^\prime)$ is the unique maximal open subset on which $s^*F$ is locally free. Similarly, $t^{-1}(U^\prime)$ is the unique maximal open subset such that the restriction of $t^*F$ is locally free. Thus, $s^{-1}(U^\prime)=t^{-1}(U^\prime)$, i.e $U^\prime\subset U$ descents to an open substack $\mathcal X^o$ of $\mathcal X$ such that $F|_{\mathcal X^o}$ is locally free.
\end{proof}

\begin{definition}[Rank of coherent sheaf]\label{def rk}
Under the hypothesis of Proposition \ref{pro rk}, we can define the rank ${\rm rk}(F)$ of $F$ to
be the rank of $F|_{\mathcal X^o}$.
\end{definition}

In order to define a notion of Gieseker stability on projective Deligne-Mumford stacks, Nironi introduced the modified Hilbert polynomial in \cite{fn}. First of all, we recall the notion of polarization on projective Deligne-Mumford stacks.
\begin{definition}[Polarization]\label{def polar}
For a Projective Deligne-Mumford stack $\mathcal{X}$, the polarization of $\mathcal X$ is a pair $(\mathcal{E},\mathcal{O}_{X}(1))$, where  $\mathcal E$ is a generating sheaf and $\mathcal{O}_X(1)$ is a very ample invertible sheaf on $X$.
\end{definition}

\begin{definition}[Modified Hilbert Polynomial]\label{def m hilbert p}
Fix a polarization $(\mathcal E,\mathcal{O}_{X}(1))$ on a projective Deligne-Mumford stack $\mathcal X$. For a coherent sheaf $F$ on $\mathcal X$, the modified Hilbert polynomial $P_F$ of $F$ is defined by
\[
P_F(m)=\mathcal{X}(\pi_\ast(F\otimes{\mathcal E^\vee})\otimes\mathcal{O}_{X}(m)), \]
where $\mathcal{X}(\pi_\ast(F\otimes{\mathcal E^\vee})\otimes\mathcal{O}_{X}(m))$ is
the Euler characteristic of $\pi_\ast(F\otimes{\mathcal E}^\vee)\otimes\mathcal{O}_{X}(m)$.
\end{definition}

\begin{remark}\label{remk hp}
In general, the modified Hilbert polynomial
\[
P_{F}(m)=\underset{i=0}{\overset{d}{\sum}}\frac{a_i(F)}{i!}\cdot m^i,
\]
where $d$ is the dimension of $F$  and $a_i(F)$ are rationals. In the special case: $F$ is a torsion free sheaf (\ref{def tfs}) on a projective algebraic orbifold $\mathcal X$ of dimension $\mathcal X$, then the coefficient $a_n(F)$ of the leading term is $\text{rk}(F)\text{rk}(\mathcal E)\text{deg}(\mathcal O_{X}(1))$, by the sake of Grothendieck-Riemann-Roch formula in \cite{bt}.
\end{remark}


\begin{definition}[Modified Slope]\label{def mslope}
Let $\mathcal X$ be an integral projective Deligne-Mumford stack over $k$. The modified Hilbert polynomial of $F$ is
$P_{F}(m)=\underset{i=0}{\overset{d}{\sum}}\frac{a_i(F)}{i!}\cdot m^i$.
The modified slope $\mu(F)$ of $F$ is
\[
\mu(F)=\frac{a_{d-1}(F)}{a_{d}(F)}.
\]
\end{definition}

Using the modified slope, we can introduce the notions of semistable (stable) torsion free sheaves.

\begin{definition}[Stability]
A torsion free sheaf $E$ is said to be semistable (resp. stable) if for all coherent subsheaves $F\subset E$ and ${\rm rk}(F)<{\rm rk}(E)$, we have
\[\mu(F)\leq\mu(E)\quad({\rm resp}.\quad\mu(F)<\mu(E)). \]
If $E$ is not semistable, $E$ is called unstable.
\end{definition}

\begin{definition}[Subbundle of torsion free sheaf]\label{def subbd}
Let $E$ be a coherent subsheaf of a torsion free sheaf $F$. If the quotient sheaf $F/E$ is also a torsion free sheaf, we say $E$ is a subbundle of $F$.
\end{definition}

Indeed, for every coherent subsheaf of a torsion free sheaf, there is a unique minimal subbundle contain it. We have the following proposition.
\begin{proposition}\label{prop subbd}
Let $\mathcal X$ be an integral projective Deligne-Mumford stack over $k$ and let $F$ be a torsion free sheaf on $\mathcal X$. For a coherent subsheaf $G$ of $F$, there is a unique coherent subsheaf $G^\prime\subseteq F$, such that
\begin{enumerate}
\item $G\subseteq G^\prime$ and ${\rm rk}(G^\prime)={\rm rk}(G)$;
\item if $F/{G^\prime}$ is not zero sheaf, then ${F}/{G^\prime}$ is a torsion free sheaf.
\end{enumerate}
\end{proposition}

\begin{proof}
We have the following two exact sequences:
\[
 \xymatrix@C=0.5cm{
   0 \ar[r] & G\ar[rr] && F \ar[rr]^{j} &&  F/G \ar[r] & 0 }, \]
\[
\xymatrix@C=0.5cm{
0 \ar[r] & T(F/G) \ar[rr]&& F/G  \ar[rr] && Q \ar[r] & 0 }, \]
where $T(F/G)$ is the maximal torsion subsheaf of $F/G$.
Then, $G^\prime = j^{-1}(T(F/G))$ and $F/G^\prime = Q$.
We have to check the uniqueness of $G^\prime$. Suppose there are two such sheaves $G_1$ and $ G_2$.
Then ${\rm rk}(G_1\cap G_2)={\rm rk}(G)$. Also, there are two exact sequences
\[
\xymatrix@C=0.5cm{
0 \ar[r] &G_1\cap G_2 \ar[rr]&& G_1 \ar[rr] && (G_1+G_2)/G_2\ar[r] & 0},
\]
\[
\xymatrix@C=0.5cm{
0 \ar[r] &G_1\cap G_2 \ar[rr]&& G_2 \ar[rr] && (G_1+G_2)/G_1\ar[r] & 0}.
\]
If $(G_1+G_2)/G_2$ and $(G_1+G_2)/G_1$ are torsion free, then
\[
{\rm rk}(G_1\cap G_2)<{\rm rk}(G_1),\quad {\rm rk}(G_1\cap G_2)<{\rm rk}(G_2).
\]
This is impossible. Hence, we have $G_1+G_2=G_1$ and
$G_1+G_2=G_2$. So, $G_1=G_2$.
\end{proof}
Following \cite{lan}, if $\mathcal X$ is a smooth algebraic orbifold, there is an explicit construction of the sheaf $G^\prime$ in Proposition \ref{prop subbd}.
\begin{proposition}\label{prop orb bundle}
For a smooth algebraic orbifold $\mathcal X$, there is an open dense substack $\mathcal X^o$ of $\mathcal X$ such that $\mathcal X^o$ is an irrducible smooth variety over $k$. Let $\xi$ be the generic point of $\mathcal X^o$ and let $\gamma:\mathcal X^o{\rightarrow}\mathcal X$ be the open immersion. Assume that $G_\xi$ and $F_\xi $ are the stalks of $\gamma^*G$ and $\gamma^*F$ at $\xi$, respectively. $G_\xi$ and $F_\xi$ can be regarded as
quasicoherent sheaves on $\mathcal X^o$. Then, we have:
\[
G^\prime=\gamma_*G_\xi\cap F.
\]
\end{proposition}
\begin{proof}
Indeed, $\gamma_*G_\xi\subseteq\gamma_*F_\xi$ and $F\subseteq\gamma_*\gamma^*F\subseteq\gamma_*F_\xi$. By $\gamma^*G\subseteq G_\xi$, we have $G\subseteq\gamma_*\gamma^*G\subseteq\gamma_*G_\xi$. Thus, $G\subseteq\gamma_*G_\xi\cap F$. Obviously, ${\rm rk}(\gamma_*G_\xi\cap F)={\rm rk}(G)$. Assume that $F/{\gamma_*G_\xi\cap F}$ is not zero. Let $\alpha:U\rightarrow\mathcal X$ be an \'etale morphism. Without loss of generality, we can assume that $U$ is an irreducible smooth affine variety ${\rm Spec}(A)$. We have the following cartesian diagram:
\begin{align*}
\xymatrix{
  U^o \ar[d]_{\alpha^o}\ar[r]^{\gamma^o} & U \ar[d]^{\alpha} \\
  \mathcal X^o \ar[r]^{\gamma} & \mathcal X   }
\end{align*}
By the flat base change theorem (Corollary A.2.2 in \cite{sb}), we only need to consider the case: $\mathcal X={\rm Spec}(A)$.
Following the Proposition 1 in \cite{lan}, if $F({\rm Spec}(A))=M$, $G({\rm Spec}(A))=N$ and $K$ is the quotient field of $A$, then $G^\prime$ is the coherent sheaf associated to the $A$-module $M\cap N\otimes_AK$.
\end{proof}
\begin{remark}\label{remk sbundle}
Under the above hypotheses, the subbundle $G^\prime$ is uniquely determined by the vector subspace $G_\xi$ of $F_\xi$ over the field of rational functions on $\mathcal X^o$.
\end{remark}

\begin{definition}[Join of sheaves]
Suppose $F_1$ and $F_2$ are two coherent subsheaves of a torsion free sheaf $F$ on a $n$-dimensional integral projective Deligne-Mumford stack $\mathcal X$ over $k$. The unique subbundle $F_1\vee F_2$ of $F$ in Proposition \ref{prop subbd} containing $F_1+F_2$, is called the join of $F_1$ and $F_2$.
\end{definition}

\begin{proposition}\label{pro inq}
If $F_1$ and $F_2$ are two subbundles of a torsion free sheaf $E$ on a $n$-dimensional integral Deligne-Mumford stack $\mathcal X$ over $k$, then the coefficients of the modified Hilbert polynomials satisfies:
\[
a_{n-1}(F_1\vee F_2)+a_{n-1}(F_1\cap F_2)\geq a_{n-1}(F_1)+ a_{n-1}(F_2).
\]
\end{proposition}

\begin{proof}

By the two exact sequences
\[
 \xymatrix@C=0.5cm{
 0 \ar[r] & F_1\cap F_2 \ar[rr]&& F_1 \ar[rr] && ( F_1+F_2)/F_2\ar[r] & 0}
\]
and
\[
 \xymatrix@C=0.5cm{
 0 \ar[r] & F_1\cap F_2 \ar[rr]&& F_2 \ar[rr] && (F_1+F_2)/F_1\ar[r] & 0}, \]
we have
\[
P_{F_1\cap F_2}+P_{(F_1+F_2)/F_1}= P_{F_2},\quad
P_{F_1\cap F_2}+P_{(F_1+F_2)/F_2}= P_{F_1}.
\]
So,
\[
a_{n-1}(F_1\cap F_2)+a_{n-1}((F_1+F_2)/F_2)= a_{n-1}(F_1)
\]
and
\[
a_{n-1}(F_1\cap F_2)+a_{n-1}((F_1+F_2)/F_1)= a_{n-1}(F_2).
\]

Also, there is an exact sequence
\[
\xymatrix@C=0.5cm{
0 \ar[r]&(F_1+F_2)/F_1\ar[rr]&&(F_1\vee F_2)/F_1 \ar[rr]&&
(F_1\vee F_2)/(F_1+F_2) \ar[r] & 0 }.
\]
Hence, we have
\[
P_{(F_1+F_2)/F_1}+P_{(F_1\vee F_2)/(F_1+F_2)}=
P_{(F_1\vee F_2)/F_1}.
\]
Therefore,
\[
a_{n-1}((F_1+F_2)/F_1)+a_{n-1}((F_1\vee F_2)/(F_1+F_2))=
a_{n-1}((F_1\vee F_2)/F_1).
\]
And also, $a_{n-1}((F_1\vee F_2)/(F_1+F_2))\geq 0$,
because $(F_1\vee F_2)/(F_1+F_2)$ is a torsion sheaf.
So,
\[
 a_{n-1}(F_1\cap F_2)+a_{n-1}((F_1\vee F_2)/F_1)\geq a_{n-1}(F_2).
\]
By the exact sequence
\[
\xymatrix@C=0.5cm{
0 \ar[r] &F_1 \ar[rr]&& F_1\vee F_2 \ar[rr] &&(F_1\vee F_2)/ F_1 \ar[r] & 0 },
\]
we have
\[
a_{n-1}(F_1\vee F_2)-a_{n-1}(F_1)= a_{n-1}((F_1\vee F_2)/F_1).
\]
Then,
\[
a_{n-1}(F_1\vee F_2)+a_{n-1}(F_1\cap F_2)\geq a_{n-1}(F_1)+ a_{n-1}(F_2).
\]
\end{proof}

As \cite{lan}, we introduce the $\beta$-invariant.
\begin{definition}
Let $E$ be a fixed torsion free sheaf on a $n$-dimensional integral projective Deligne-Mumford stack $\mathcal X$ over $k$. For every torsion free sheaf $F$ on $\mathcal X$, we can define
the $\beta$-invariant as
\[
\beta(F)=a_n(E)a_{n-1}(F)-a_{n-1}(E)a_n(F).
\]
\end{definition}

\begin{remark}
By Proposition \ref{prop subbd}, if every proper subbundle $F\subset E$ satisfies $\beta(F)\leq 0$, then $E$ is semistable.
\end{remark}

\begin{proposition}\label{prop pro beta}
Let $\mathcal X$ be an integral projective Deligne-Mumford stack over $k$.
\begin{enumerate}
\item If $F_1$ and $F_2$ are two subbundles of $E$ on $\mathcal X$, then
\[
\beta(F_1)+\beta(F_2)\leq \beta(F_1\vee F_2) + \beta(F_1\cap F_2),
\]
with equality if and only if the codimension of the sheaf $(F_1\vee F_2)/(F_1+F_2)$ $\geq 2$.
\item If $\xymatrix@C=0.5cm{ 0 \ar[r] & F \ar[r] & G \ar[r]& K \ar[r] & 0 }$ is exact sequence
 of torsion free sheaves on $\mathcal X$, then
\[
\beta(F)+\beta(K)=\beta(G). \]
\end{enumerate}
 \end{proposition}

\begin{proof}
For the first statement, we have
\[
 \beta(F_1)+\beta(F_2)=a_n(E)a_{n-1}(F_1)-a_{n-1}(E)a_n(F_1)+
a_n(E)a_{n-1}(F_2)-a_{n-1}(E)a_n(F_2)=
\]
\[a_n(E)\big(a_{n-1}(F_1)+a_{n-1}(F_2)\big)-a_{n-1}(E)\big(a_n(F_1)+a_n(F_2)\big)\leq \]
\[
a_n(E)\big(a_{n-1}(F_1\vee F_2)+a_{n-1}(F_1\cap F_2)\big)-a_{n-1}(E)\big(a_n(F_1)+a_n(F_2)\big).
\]
By the exact sequence $\xymatrix@C=0.5cm{0 \ar[r] & F_1\cap F_2 \ar[rr]&& F_1\oplus F_2
\ar[rr] && F_1+F_2\ar[r] & 0 }$, we have
\[
P_{F_1\oplus F_2}=P_{F_1\cap F_2} + P_{F_1+F_2}.
\]
So,
\[                                                                                                                            a_{n}(F_1+F_2)+a_{n}(F_1\cap F_2)= a_{n}(F_1)+ a_{n}(F_2).
\]
Also, $a_n(F_1+F_2)=a_n(F_1\vee F_2)$. Then,
\[
\beta(F_1)+\beta(F_2)\leq \beta(F_1\vee F_2) + \beta(F_1\cap F_2)
\]
The second statement is obvious.
\end{proof}

Following \cite{lan}, we consider the set $\Gamma(E)$ of proper subbundles of $E$, which have the following Property:
\[
\Gamma(E)=\{F: \text{$F$ is a proper subbundle of $E$ such that for every subsheaf $G\subset F$, $\beta(G)<\beta(F)$}\}.
\]

\begin{remark}
The set $\Gamma(E)$ is nonempty. In fact, the zero sheaf is in $\Gamma(E)$. In addition, if $E$ is semistable, there is only one element in the set $\Gamma(E)$, i.e the zero sheaf.
\end{remark}

\begin{proposition}\label{prop beta bd}
Let $F$ be a maximal element of $\Gamma(E)$. For every subbundle $G\supseteq F$, we have
$\beta(G)\leq\beta(F)$.
\end{proposition}
\begin{proof}
Suppose $\beta(G)>\beta(F)$. Let $H\subset G$ be the minimal subbundle such that
$\beta(H)>\beta(F)$ and $F\subseteq H$. For every proper subbundle $I$ of $H$ and
$F\nsubseteq I$, we have
\[
 \beta(I\vee F)-\beta(I)\geq \beta(F)-\beta(F\cap I)>0. \]
By the minimality of $H$, $\beta(H)\geq\beta(I\vee F)$. So, $\beta(H)>\beta(I)$.
Therefore, $H\in\Gamma(E)$. Contradiction!
\end{proof}

\begin{corollary}\label{cor beta}
 There is unique maximal subbundle $F\in\Gamma(E)$. Also, for every subbundle $B\subseteq E$, $\beta(B)\leq\beta(F)$ with equality only if $B\supseteq F$.
\end{corollary}
\begin{proof}
If there are two maximal subbundles $F_1$ and $F_2$ in $\Gamma(E)$, then
\[
\beta(F_1\vee F_2)-\beta(F_1)\geq\beta(F_2)-\beta(F_1\cap F_2).
\]
By Proposition \ref{prop beta bd}, $\beta(F_1\vee F_2)\leq\beta(F_1)$. Thus, $\beta(F_2)\leq\beta(F_1\cap F_2)$. On the other hand, $\beta(F_1\cap F_2)\leq\beta(F_2)$. Then, $\beta(F_1\cap F_2)=\beta(F_2)$. So, $F_1\cap F_2=F_2$. Similarly, $F_1\cap F_2=F_1$. Then, $F_2=F_1$. Hence, there is a unique maximal subbundle $F\in\Gamma(E)$. By $\beta(F\vee B)-\beta(B)\geq\beta(F)-\beta(F\cap B)\geq 0$ and $\beta(F)\geq\beta(F\vee B)$, we get $\beta(F)\geq\beta(B)\text{ with equality only if $B\supseteq F$}.$
\end{proof}

\begin{remark}
In the above corollary, the $\beta(F)$ is the maximum value of $\beta$-invariant for subbundles in $E$. Also,
${\rm Hom}_{\mathcal O_X}(F,E/F)$=0.
\end{remark}
At the end of this section, we show that the torsion free semistable sheaves is stable under the extension of the base field $k$($k$ is not necessarily algebraic closed).
\begin{proposition}\label{pro ex ss}
Let $k^\prime$ be an extension field of $k$. We have the following cartesian diagram:
\[
\xymatrix{
  \mathcal X\times{\rm{Spec}}(k^\prime) \ar[d]_{p_2} \ar[r]^{\qquad p_1} &\mathcal X \ar[d] \\
  {\rm{Spec}}(k^\prime) \ar[r] & {\rm{Spec}}(k)}
\]
Assume that the field $k$ is infinite when $k^\prime/k$ is not algebraic. Then $E^\prime=p_1^*E$ is semistable if and only if $E$ is semistable.
\end{proposition}
\begin{proof}
The proof can be proved as Proposition 3 \cite{lan}, or can be found in \cite{fn}.
\end{proof}
\section{The Main Results}
From now on, $\mathcal X$ is a $n$-dimensional smooth algebraic orbifold with a fixed polarization $(\mathcal E,\mathcal O_{\mathcal X}(1))$ over $k$. Let $R\supseteq k$ be a discrete valuation ring with maximal ideal $m=(\pi)$ and residue field $k$. $K$ is the quotient field of $R$. Consider the following cartesian diagram:
\begin{equation*}
\xymatrix{
\mathcal X_K \ar@{^(->}[r]^i \ar[d] & \mathcal X_R \ar[d] & \mathcal X_k \ar[d] \ar@{_(->}[l]_{j}\\
     \text{Spec}(K) \ar@{^(->}[r] & \text{Spec}(R) & \text{Spec}(k) \ar@{_(->}[l] }
\end{equation*}
where $\mathcal X_R=\mathcal X\times \text{Spec}(R)$, $\mathcal X_K=\mathcal X\times \text{Spec}(K)$ and
$\mathcal X_k=\mathcal X\times\text{Spec}(k)=\mathcal X$. $i:\mathcal X_K\rightarrow\mathcal X_R$ is the natural open immersion and $j:\mathcal X_k\rightarrow\mathcal X_R$ is the natural closed immersion. Our goal is to prove the following result:
\begin{theorem}\label{thm main 1}
Assume that $E_K$ is a torsion free sheaf on $\mathcal X_K$. Then
\begin{enumerate}
\item If $E_1$ and $E_2$ are two coherent subsheaves of $i_*E_K$ on $\mathcal X_R$ such that $i^*E_1=i^*E_2=E_K$ and $j^*E_1$, $j^*E_2$ are semistable torsion free sheaves on $\mathcal X_k$, at least one of which is stable, then there is an integer $p$ such that $E_1=\pi^p E_2$.
\item If $E_K$ is semistable, then there exists a coherent subsheaf $E\subseteqq {i_*}E_K$ such that $i^*E=E_K$ and $j^*E$ is torsion free and semistable on $\mathcal X_k$.
\end{enumerate}
\end{theorem}

We first state a lemma, which corresponds to Proposition 5 in \cite{lan}.
\begin{lemma}\label{lemm shp}
If $E_1$ and $E_2$ are two torsion free sheaves on $\mathcal X_{R}$ such that $i^*E_1=i^*E_2$, then the modified Hilbert polynomials $P_{j^*E_1}(m)=P_{j^*E_2}(m)$. In particular, $a_{n-1}(j^*E_1)=a_{n-1}(j^*E_2)$.
\end{lemma}
\begin{proof}
Since the field $k$ is algebraically closed and $R$ is a regular local ring, $\mathcal X_R$ is integral and smooth over $\text{Spec}(R)$. Then, the torsion free sheaf on $\mathcal X_R$ is flat over $\text{Spec}(R)$, since the torsion free modules over valuation rings are flat. By the Lemma 3.16 in \cite{fn}, $P_{j^*E_1}(m)=P_{j^*E_2}(m)$.
\end{proof}
$\mathcal X$ has an open dense substack $\mathcal X^o$ such that it is an irreducible smooth variety over $k$. Let $\gamma:\mathcal X^o\rightarrow \mathcal X$ be the corresponding open immersion. $\mathcal X_K^o=\mathcal X^o\times\text{Spec}(K)$ and
$\mathcal X_k^o=\mathcal X\times\text{Spec}(k)$ are also irreducible and smooth. Let $\Xi$ be the generic point of $\mathcal X_K^o$
and $\xi$ be the generic point of $\mathcal X_k^o$. And, we have the following cartesian diagram:
\[
 \xymatrix{
   \mathcal X_K^o \ar[d]  \ar[r] & \mathcal X_R^o \ar[d]  & \mathcal X_k^o \ar[l] \ar[d]\\
   \mathcal X_K \ar[r] & \mathcal X_R   & \mathcal X_k \ar[l]  }
\]
Let $E_K$ be a torsion free sheaf of rank $r$ on $\mathcal X_K$. Since $\mathcal X^o$ is an integral scheme, the stalk $(E_K)_{\Xi}$ of $(E_K)|_{\mathcal X^o_K}$ at $\Xi$ is a free $\mathcal O_{\Xi}$ module. Denote the stalks of $\mathcal O_{\mathcal X^o_R}$ at $\Xi$ and $\xi$ by $\mathcal O_{\Xi}$ and $\mathcal O_{\xi}$, respectively.

\begin{lemma}\label{lemm ext 1}
Suppose $M\subset(E_K)_{\Xi}$ is a free rank $r$ $\mathcal O_\xi$-submodule of $(E_K)_{\Xi}$. Then there exists a unique torsion free sheaf $E\subseteq i_*E_K$ on $\mathcal X_R$ such that $i^*E=E_K$, $E_\xi=M$, and $j^*E$ is a torsion free sheaf on $\mathcal X_k$.
\end{lemma}

\begin{proof}
As above, $\mathcal O_\xi$ is the stalk of $\mathcal O_{\mathcal X_R^o}$ at the generic $\xi$ of $\mathcal X_k^o$ in $\mathcal X_R^o$. Then, there is a natural morphism $\beta_1:\rm Spec (\mathcal O_\xi)\rightarrow \mathcal X_R^o$. Besides, $\Xi$ is the generic point of $\mathcal X_K^o$. So, $\Xi$ is also the generic point of $\mathcal X_R^o$. So, there are two natural morphisms $\alpha:\Xi\rightarrow \rm{Spec}(\mathcal O_\xi)$ and $\beta_2:\Xi\rightarrow\mathcal X_K^o$. Let $i^o:\mathcal X_K^o\rightarrow \mathcal X_R^o$ be the open immersion obtained through base change from the open immersion ${\rm Spec}(K)\hookrightarrow {\rm Spec}(R)$. And also, they form the following cartesian diagrams:
\begin{equation*}\label{diag 1}
\begin{split}
\xymatrix{
  \text{Spec}(\mathcal O_{\xi})\ar[r]^{\quad\beta_1} & \mathcal X_R^o  \ar[r]^{\gamma_R} & \mathcal X_R \\
  \Xi \ar[r]^{\beta_2} \ar[u]_{\alpha} & \mathcal X_K^o \ar[u]_{i^o} \ar[r]^{\gamma_K} &
  \mathcal X_K \ar[u]_{i} }
\end{split}\tag{\rm{A}}
\end{equation*}
Denote the torsion free sheaf on $\rm Spec(\mathcal O_\xi)$ corresponding to the modules $M$ by $\mathcal M$. Similarly, $\mathcal N$ is the free sheaf on $\Xi$ corresponding to $(E_K)_{\Xi}$.
\begin{center}\label{claim 1}
\textbf{Claim}: $E=i_*E_K\cap(\gamma_R\circ\beta_1)_*\mathcal M$ satisfies the conditions in the conclusion of Lemma \ref{lemm ext 1}.
\end{center}
\textbf{First step}: we need to explain the intersection of $i_*{E_K}$ and $(\gamma_R\circ\beta_1)_*\mathcal M$ in $(\gamma_R\circ\beta_1\circ\alpha)_*\mathcal N$.\\
By the inclusion $M\subseteq(E_K)_{\Xi}$, we have the inclusion:
\begin{equation*}\label{inclu a}
{(\gamma_R\circ\beta_1)}_*\mathcal M \subseteq (\gamma_R\circ\beta_1)_*(\alpha_*\mathcal N).\tag{a}
\end{equation*}
In addition, $\Xi$ is the generic point of $\mathcal X_K^o$, there is another inclusion:
\begin{equation*}\label{inclu b}
(i\circ{\gamma_K})_*({\gamma_K}^*E_K)\subseteq (i\circ\gamma_K)_*({\beta_2}_*\mathcal N).\tag{b}
\end{equation*}
By the diagram (\ref{diag 1}), we get ${(i\circ\gamma_K\circ\beta_2)}_*\mathcal N ={{(\gamma_R\circ\beta_1\circ\alpha)}_*}\mathcal N$. In the following, we show that the morphism:
\begin{equation}\label{equ sec 1}
E_K\longrightarrow {\gamma_K}_*({\gamma_K}^*E_K)
\end{equation}
obtained by adjunction formula is injective. Indeed, the coarse moduli space of $\mathcal X_K$ is $X_K=X\times {\rm Spec}(K)$ and $X_K$ is irreducible. By the Remark \ref{remk irreducible}, $\mathcal X_K$ is irreducible. Also, $\mathcal X_K$ is reduced. Then, $\mathcal X_K$ is integral. For every \'etale morphism $f:U\rightarrow\mathcal X$ from an irreducible smooth variety $U$ over $k$ to $\mathcal X$, we have the cartesian diagram:
\begin{equation*}\label{diag 2}
\begin{split}
\xymatrix{
  U^o_K \ar[d]_{f_K^o} \ar[r]^{\gamma_K} & U_K \ar[d]^{f_K} \\
  \mathcal X^o_K  \ar[r]^{\gamma^\prime_K}   & \mathcal X_K }
\end{split}\tag{\rm{B}}
\end{equation*}
where $U_K=U\times{\rm Spec}(K)$. Pulling back the homomorphism (\ref{equ sec 1}) to $U_K$, we get
\begin{equation}\label{equ sec 2}
f_K^*E_K\longrightarrow f_K^*{\gamma_K}_*({\gamma_K}^*E_K).
\end{equation}
By the flat base change theorem of stacky version ( Corollary \rm{A}.2.2 in \cite{sb} and \rm{A}.3.4 in \cite{sb1}), we have
\begin{equation}\label{equ sec 3}
f_K^*{\gamma_K}_*({\gamma_K}^*E_K)={\gamma_K^\prime}_*{f^o_K}^*({\gamma_K}^*E_K).
\end{equation}
On the other hand, ${\gamma_K^\prime}^*f_K^*E_K={f^o_K}^*{\gamma_K}^*E_K$. Then, the homomorphism (\ref{equ sec 2}) is
\begin{equation}\label{equ sec 4}
f_K^*E_K\longrightarrow{\gamma_K^\prime}_*{\gamma_K^\prime}^*f_K^*E_K .
\end{equation}
Because $U$ is integral and $f_K^*E_K$ is torsion free, the homomorphism (\ref{equ sec 4}) is injective. Thus, the homomorphism (\ref{equ sec 1}) is injective. So, $i_*{E_K}\longrightarrow i_*{\gamma_K}_*{\gamma_K}^*E_K$ is injective.
Hence, by (\ref{inclu a}), (\ref{inclu b}) and the diagram (\ref{diag 1}), we have the following two short exact sequences with the same middle terms:
\begin{equation}
\xymatrix{
&  &  0 \ar[d]                                  \\
&  &{(\gamma_R\circ\beta_1)}_*\mathcal M \ar[d] \\
0 \ar[r] & i_*E_K \ar[r] & {{(\gamma_R\circ\beta_1\circ\alpha)}_*}\mathcal N}
\end{equation}
Thus, $E =i_*E_K\cap{(\gamma_1\circ\beta_1)}_*\mathcal M$ is a quasicoherent sheaf on $\mathcal X_R$. We accomplished the first part of the proof.\\
\textbf{Second step}: We have to check the sheaf $E$ which we have defined is a torsion free coherent sheaf. We only need to check this locally in the \'etale topology.
Suppose $\theta:{\rm{Spec}}(A)\rightarrow \mathcal X$ is an \'etale morphism and ${\rm{Spec}}(A)$ is a smooth irreducible variety over $k$. We have the cartesian diagram
\[
\xymatrix{
  \mathcal X^o \ar[r] & \mathcal X \\
   V \ar[r]\ar[u]^{\phi} & \text{Spec}(A)\ar[u]^{\theta} }
\]
Since $\phi$ is an \'etale morphism of finite type between irreducible smooth varieties, $\phi$ is generically finite dominant map i.e $\phi^{-1}(\xi)$ is a finite set. By exercise 3.7 in page 91 of \cite{Ha},
there is an open dense subset $i_W:W\rightarrow\mathcal X^o$ such that the morphism
$\phi^{\prime}:\phi^{-1}(W)\rightarrow W$ is finite and \[
\xymatrix{
   W \ar[r]^{i_W} & \mathcal X^o \ar[r] & \mathcal X \\
 W^{\prime}=\phi^{-1}(W) \ar[r] \ar[u]^{\phi^\prime}&  V \ar[r] \ar[u]^{\phi} & \text{Spec}(A)\ar[u]^{\theta}  }
\]
is a cartesian diagram. Denote $\phi^{-1}(W)$ by $W^{\prime}$. By the base change, we have
\begin{equation}\label{diag 3}
\begin{split}
\xymatrix{
\Xi\ar[r]^{\alpha \qquad}& \text{Spec}(\mathcal O_{\xi}) \ar[r]^{\beta_3} &W\times\text{Spec}(R)\ar[r]^{\qquad i_{W,R}}&\mathcal X_R^o\ar[r]^{\gamma_1}&\mathcal X_R\\
\Xi_1 \ar[u]^{\phi_{\Xi}} \ar[r]^{\alpha^\prime\qquad}&\text{Spec}(\mathcal O_{V_R,\xi^\prime})\ar[r]^{\beta_3^\prime}
\ar[u]^{\phi_{\xi}} & W^{\prime}\times\text{Spec}(R) \ar[r]^{i_{W^\prime,R}} \ar[u]^{\phi^{\prime}_R} &V\times\text{Spec}(R)\ar[r]^{\gamma_1^\prime}\ar[u]^{\phi_R} & \text{Spec}(A\otimes_kR)\ar[u]^{\theta_R}
}
\end{split}\tag{C}
\end{equation}
where $\Xi^{\prime}$ is the generic point of $V_R=V\times{\rm Spec}R$ and $\xi^\prime$ is the generic point of the close subscheme $W^{\prime}\times{\rm Spec}(k)\hookrightarrow W^{\prime}\times{\rm Spec}(R)$.
The first square and the second square are cartesian. Indeed, we may assume $W=\text{Spec}(B)$ and $W^{\prime}=\text{Spec}(C)$. Then ${\phi^\prime}^\sharp:B\rightarrow C$ is an injective finite map. $\xi$ and $\xi^{\prime}$ are the prime ideals $B\otimes_{k}(\pi)$ and $C\otimes_{k}(\pi)$ respectively. Denote the quotient fields of $B$ and $C$ by $K_B$ and $K_C$ respectively. Since the field $k$ is algebraically closed, it follows that $\mathcal O_{\xi}= K_B\otimes_{k}R$ and $\mathcal O_{V_R,\xi^\prime}=K_C\otimes_{k}R$. Then
\[
\mathcal O_{\xi}\otimes_{B\otimes_{k}R}(C\otimes_{k}R)=
(K_B\otimes_kR)\otimes_{B\otimes_kR}(C\otimes_kR)=
(K_B\otimes_{B}(B\otimes_{k}R))\otimes_{B\otimes_{k}R}(C\otimes_{k}R)=
\]
\[
K_B\otimes_{B}(C\otimes_{k}R)=(K_B\otimes_{B}C)\otimes_{k}R=
K_C\otimes_{k}R=\mathcal O_{V_R,\xi^\prime},
\]
where $K_C=K_B\otimes_{B}C$ ($C$ is integral over $B$). Thus, the second square is cartesian. So, the morphism $\phi_\xi$ is finite. Then the first square is cartesian. By the flat base change formula of stacky version and cartesian diagram
\begin{equation*}\label{diag 4}
\begin{split}
\xymatrix{
 \Xi \ar[r] &W\times\text{Spec}(K) \ar[r]^{\qquad i_{U,K}} &\mathcal X^o_K \ar[r] &\mathcal X_K \ar[r]^{i} & \mathcal X_R \\
\Xi^{\prime} \ar[r] \ar[u]^{\phi_{\Xi}}&W^\prime\times\text{Spec}(K) \ar[r] \ar[u]^{\phi^\prime_K}
&V\times\text{Spec}(K)\ar[r] \ar[u]^{\phi_K}&
\text{Spec}(A\otimes_{k}K)\ar[r]^{i^\prime} \ar[u]^{\theta_K} &
\text{Spec}(A\otimes_{k}R)\ar[u]^{\theta_R}
}
\end{split}\tag{D}
\end{equation*}
By the last square in diagram (\ref{diag 4}), we have the equation:
\begin{equation}\label{equ sec 5}
 \theta_R^*E=\theta_R^*\big(i_*E_K\cap{(\gamma_1\circ\beta_1)}_*\mathcal M\big).
\end{equation}
From the last three square in diagram (\ref{diag 3}) , we get the equation:
\begin{equation}\label{equ sec 6}
\theta_R^*i_*E_K\cap\theta_R^*\big({(\gamma_1\circ\beta_1)}_*\mathcal M\big)=
i_*^\prime\theta_K^*E_K\cap {(\gamma_1^\prime\circ i_{W^\prime,R}\circ\beta_3^\prime)}_*\phi_{\xi}^*\mathcal M.
\end{equation}
Let $\phi_{\xi}^*\mathcal M=\mathcal M^\prime$, $\theta_K^*E_K=E^\prime$ and $\phi_{\Xi}^{*}\mathcal N =
\mathcal N^\prime$. Then,
\begin{enumerate}
\item $E^\prime$ is a torsion free sheaf of rank $r$ and $E^\prime|_{\Xi}=\mathcal N^\prime$;
\item ${\alpha_1^{\prime}}^*\mathcal M^\prime$ = $\mathcal N^\prime$; \item $\mathcal M^\prime$ and $\mathcal N^\prime$ are free sheaves of rank $r$.
\end{enumerate}
Therefore, we only consider the case: $\mathcal X={\rm{Spec}}(A)$, where ${\rm{Spec}}(A)$ is an irreducible smooth affine varieties over $k$. In this case, all the properties of $E$ can be checked through commutative algebra, just as the Proposition 6 of \cite{lan}.
\end{proof}

\begin{remark}\label{rmk of ext}
Assume that $M_1$ and $M_2$ are two free rank $r$ $\mathcal O_{\xi}$ submodules of $(E_K)_{\Xi}$. Denote the corresponding coherent sheaves in Lemma \ref{lemm ext 1} by $E_1$ and $E_2$, respectively. If $M_1\subseteq M_2$, by the proof of Lemma \ref{lemm ext 1}, we have $E_1\subseteq E_2$.
\end{remark}
In the following, we show the first part of Theorem \ref{thm main 1} as \cite{lan}.
\begin{proof}[\textbf{The Proof of the first part in Theorem \ref{thm main 1}}]
Suppose $E_1$ and $E_2$ are two coherent subsheaves of $i_*E_K$ such that $i^*E_1=i^*E_2=E_K$ and $j^*E_1$, $j^*E_2$ are  torsion free semistable sheaves on $\mathcal X_k$, at least one of which is stable. Since $\mathcal O_{\xi}$ is a principal ideal domain, $E_{1,\xi}$ and $E_{2,\xi}$ are free $\mathcal O_{\xi}$ modules of rank $r$. Also,
$E_{1,\xi}\otimes_{\mathcal O_{\xi}}\mathcal O_{\Xi}=E_{2,\xi}\otimes_{\mathcal O_{\xi}}\mathcal O_{\Xi}=(E_K)_{\Xi}$. By the elementary divisor theorem (Theorem 7.8 in \cite{sl}), there is a basis $\{e_1,\ldots,e_r\}$ of $E_{1,\xi}$ over $\mathcal O_{\xi}$ such that $\{\pi^{q_1}e_1,\ldots,\pi^{q_r}e_r\}$ is a basis of $E_{2,\xi}$. Since we are trying to prove that $E_1=\pi^pE_2$ for some $p$, we may multiply $E_{2,\xi}$ by $\pi^m$ for some integer $m$, so that all the $q_i$ are nonnegative, and at least one of the $q_i=0$. If all the $q_i=0$ , we are done; hence we may also assume that some $q_i$ is postive. By $E_{2,\xi}\subseteq E_{1,\xi}$ and the Remark \ref{rmk of ext}, we have $E_2\subseteq E_1$. This inclusion induces a homomorphism $\alpha:j^*E_2\rightarrow j^*E_1$ on $\mathcal X_k$. Also, $\text{rk}(j^*E_1)=\text{rk}(j^*E_2)$ and
$a_{n-1}(j^*E_1)=a_{n-1}(j^*E_2)$, for the sake of the Lemma \ref{lemm shp}. Hence, $j^*E_1$ and $j^*E_2$ have the same modified slope. By the construction of $\alpha$, the map $\alpha$ is not zero and not isomorphism in codimension one. Therefore, we have $E_1=\pi^pE_2$, for some integer $p$.
\end{proof}
We state a Lemma about the torsion free modules on a discrete valuation ring.
\begin{lemma}\label{lemm tor dvr}
Suppose $M$ is a finitely generated torsion free module on a discrete valuation ring. Then $M$ is a free module of finite rank.
\end{lemma}
On analogy with \cite{lan}, we introduce Bruhat-Tits complex of the $E_K$.
Assume that $\mathfrak M$ is the set of all free rank $r$ $\mathcal O_\xi$ submodules of $(E_K)_\Xi$. For every $M\in\mathfrak M$, there is a unique torsion free sheaf $E_R$ on $\mathcal X_R$, which is the extension of $E_K$, for the sake of Lemma \ref{lemm ext 1}. An equivalence relation $\sim$ is defined in $\mathfrak M$ by
\begin{center}\label{claim 1 }
For $M,M^\prime\in\mathfrak M$, then $M\sim M^\prime$ if and only if $M=\pi^pM^\prime$, for some $p\in\mathbb Z$. \qquad(E)
\end{center}
Let $\mathfrak Q$ be the set of equivalence classes in $\mathfrak M$. Obviously, every equivalence class in $\mathfrak Q$, defines an extension of $E_K$ to coherent sheaf on $\mathcal X_R$, modulo isomorphism. We now define the structure of an $r$-dimensional simplicial complex on $\mathfrak Q$, which we will call the Bruhat-Tits complex. The dimension of $\mathfrak Q$ will be less than or equal to $r$. Two equivalence classes $[M]$ and $[M^\prime]$ in $\mathfrak Q$ are said to be adjacent if $M$ has a direct decomposition $M=N\oplus P$ such that $M^\prime=N+\pi M$. Since $\mathcal O_{\xi}$ is a discrete valuation ring, $M$ has a basis $\{e_1,e_2,\ldots,e_r\}$ over $\mathcal O_\xi$ such that $\{e_1,\ldots,e_s\}$ and $\{e_{s+1},\ldots,e_r\}$, are bases of $N$ and $P$, respectively, by the sake of Lemma \ref{lemm tor dvr}. So, $\{e_1,\ldots,e_s,\pi e_{s+1},\ldots,\pi e_r\}$ is a basis of $M^\prime$ over $\mathcal O_\xi$. Then, $M$ is adjacent to $M^\prime$ if and only if there is a basis $\{e_1,\ldots,e_r\}$ of $M$ such that $\{e_1,\ldots,e_s,\pi e_{s+1},\ldots,\pi e_r\}$ is a basis of $M^\prime$. A chain $0\subset N_1\subset N_2\subset\cdots\subset N_i\subset M$ of submodules such that each $N_i$ is a direct factor of $M$ and $M_i=N_i+\pi M$, then the $i+1$ mutually adjacent vertices  $[M],[M_1],\ldots,[M_i]$ are said to form a $i$-simplex in $\mathfrak Q$. In other words, the vertices $[M],[M_1],\ldots,[M_i]$ are said to form a $r$-simplex in $\mathfrak Q$ if there is a basis $\{e_1,e_2,\ldots,e_r\}$ of $M$ such that $N_k=(e_1,\ldots,e_{s_k})$ and $M_k=(e_1,\ldots,e_{s_k}, \pi e_{s_k+1},\ldots, \pi e_r)$, for $1\leq k\leq i$. From the above argument, it is clear that the proof of the part $2$ in the Theorem \ref{thm main 1} is equivalent to find a vertex $[E_\xi]$ of $\mathfrak Q$ such that the reduction $E_k$ of the corresponding extension $E_R$ is semistable. Start with any vertex $[E_\xi]$ in $\mathfrak Q$. We have the following Proposition, which is the orbifold vertion of Proposition 7 in \cite{lan}.

\begin{proposition}\label{prop b t c}
Assume that $[E_\xi]$ is a vertex in $\mathfrak Q$ and $E_k$ is the corresponding sheaf on $\mathcal X_k$. Then, there is a natural one-to-one correspondence between edges in $\mathfrak Q$ at $[E_\xi]$ and proper subbundles of $E_k$.
Furthermore, if $F\subset E_k$ is a subbundle corresponds to the edge $[E_\xi]-[E_\xi^\prime]$, and if $Q^\prime\subset E_k^\prime$ is the subbundle corresponds to the edge $[E_\xi^\prime]-[E_\xi]$ at $[E_\xi^\prime]$, then there are a homomorphism $E_k\rightarrow E_k^\prime$ with kernel $F$ and image $Q^\prime$, and a homomorphism
$E^\prime\rightarrow E_k$ with kernel $Q^\prime$ and image $F$.
\end{proposition}

\begin{proof}
First, let $E_\xi=(e_1,\ldots,e_r)$ be a representative of the given vertex $[E_\xi]$ and let $E_\xi^\prime=(e_1,\ldots,e_s,\pi e_{s+1},\ldots, \pi e_r)$ be a representative of an adjacent vertex. By the Remark \ref{rmk of ext}, we have a natural inclusion of the corresponding extensions $E_R^\prime$ into $E_R$. If $\widehat{E_\xi}$ and $\widehat{E^\prime_\xi}$ are the coherent sheaves on ${\rm Spec}(\mathcal O_\xi)$, defined by $E_\xi$ and $E_\xi^\prime$, respectively. In the proof of Lemma \ref{lemm ext 1}, we have showed that $E_R=i_*E_K\cap(\gamma_R\circ\beta_1)_*\widehat{E_\xi}$ and  $E_R^\prime=i_*E_K\cap(\gamma_R\circ\beta_1)_*\widehat{E_\xi^\prime}$, where the morphisms $\gamma_R$ and $\beta_1$ are the same as in diagram (\ref{diag 1}). Let $Q_\xi$ be the cokernel of inclusion $E_\xi^\prime\hookrightarrow E_\xi$ and let $\xymatrix@C=0.5cm{
0 \ar[r] & \widehat{E_\xi}\ar[r] & \widehat{E_\xi^\prime} \ar[r] & \widehat{Q_\xi} \ar[r] & 0 }$ be the associated exact sequence of coherent sheaves on ${\rm Spec}(\mathcal O_\xi)$. By the cartesian diagram (\ref{diag 3}), the following sequence:
\[
\xymatrix@C=0.5cm{
  0 \ar[r] &  (\gamma_R\circ\beta_1)_*\widehat{E_\xi^\prime}\ar[rr] && (\gamma_R\circ\beta_1)_*\widehat{E_\xi}\ar[rr] && (\gamma_R\circ\beta_1)_*\widehat{Q_\xi}\ar[r] & 0 }
\]
is exact. Thus, the cokernel $Q$ of $E_R^\prime\hookrightarrow E_R$ admits an injection $Q\hookrightarrow (\gamma_R\circ\beta_1)_*\widehat{Q_\xi}$. So, $Q$ is a coherent $\mathcal O_{\mathcal X_k}$ module. Restricting to $\mathcal X_k$, we get right exact sequence:
\[
\xymatrix@C=0.5cm{
   E_k^\prime \ar[rr] &&  E_k  \ar[rr] && Q \ar[r] & 0 }.
\]
And also, $Q$ is torsion free on $\mathcal X_k$. Indeed, as the proof of Lemma \ref{lemm ext 1}, we only need to check this, when $\mathcal X$ is an irreducible smooth affine variety. Assume that $\mathcal X={\rm Spec}(A)$. Then, $(\gamma_R\circ\beta_1)_*\widehat{Q_\xi}$ is isomorphic to the quasicoherent sheaf that associated to the direct sum of $(r-s)$-copies $K_A$, where $K_A$ is the quotient field of $A$. Thus, the image $F=\text{Im}(E_k^\prime\rightarrow E_k)$ is a subbundle of $E_k$, with an exact sequence:
\[
\xymatrix@C=0.5cm{
  0 \ar[r] & F \ar[rr] && E_k \ar[rr] && Q \ar[r] & 0 }.
\]
Now, we have construct a subbundle $F$ of $E_k$, from an edge at $[E_\xi]$.\\
\hspace*{2em} Conversely, if $F$ is a subbundle of $E_k$ and $Q=E_k/F$, then we have an exact sequence of torsion free sheaves:
\[
\xymatrix@C=0.5cm{
  0 \ar[r] & F \ar[rr] && E_k \ar[rr] && Q \ar[r] & 0 }.
\]
On the other hand, there is a natural surjective homomorphism $E_R\rightarrow E_k$. Composing this morphism with
the last morphism in the above exact sequence, we get a surjective homomorphism $E_R\rightarrow Q$ of coherent sheaves and an exact sequence:
\begin{equation}\label{exact seq 1}
\xymatrix@C=0.5cm{
  0 \ar[r] & E^\prime_R \ar[rr] && E_R \ar[rr] && Q \ar[r] & 0 }.
\end{equation}
We have to show that the above two procedures are invertible to each other. In fact, by the exact sequence (\ref{exact seq 1}), we have:
\[
\xymatrix{
  0 \ar[r] & ({E^\prime_R})_\xi  \ar[rr] && ({E_R})_\xi \ar[dr] \ar[rr] && Q_\xi \ar[r] & 0 \\
           &                             &&                            &({E_k})_\xi \ar[ur]&
}
\]
Suppose that $(E_k)_\xi$ is generated by $\{\overline e_1,\ldots,\overline e_r\}$ and $F_\xi$ is
generated by $\{\overline e_1,\ldots,\overline e_s\}$. Also, $\{\overline e_1,\ldots,\overline e_r\}$ lifts to a basis
$\{e_1,\ldots,e_r\}$ of $(E_R)_\xi$. Then, $({E^\prime_R})_\xi$ is generated by $\{e_1,\ldots,e_s, \pi e_{s+1},\ldots, \pi e_r\}$.
And, $({E^\prime_R})_\xi$ represents a vertex of $\mathfrak Q$, which is adjacent to $[E_\xi]$. Pulling back the exact sequence
\[
\xymatrix@C=0.5cm{
  0 \ar[r] & E^\prime_R \ar[rr] && E_R \ar[rr] &&  Q \ar[r] & 0 }
\]
to $\mathcal X_k$, we get

\[
\xymatrix@C=0.5cm{
  0 \ar[r] & Q^\prime\ar[r]& E^\prime_k \ar[r] & E_k \ar[r] &  Q \ar[r] & 0 },
\]
where $Q^\prime=\text{Tor}_1^{\mathcal O_{\mathcal X_R}}(Q,\mathcal O_{\mathcal X_k})$. Tensoring the exact sequence $\xymatrix@C=0.5cm{0 \ar[r] &
\mathcal O_{\mathcal X_R} \ar[r]^{\pi} &\mathcal O_{\mathcal X_R}\ar[r] & \mathcal O_{\mathcal X_k} \ar[r] & 0 }$ with $Q$, we have
\[
\xymatrix@C=0.5cm{
  0 \ar[r] & Q^\prime\ar[r]& Q \ar[r]^{\pi} & Q \ar[r]^{\text{id}} & Q \ar[r] & 0 },
\]
whence $Q^\prime\cong Q$. Thus, we get two exact sequences
\[
\xymatrix@C=0.5cm{
  0 \ar[r] & F \ar[rr] && E_k \ar[rr] && Q \ar[r] & 0 };
\]
\[
\xymatrix@C=0.5cm{
  0 \ar[r] & Q   \ar[rr] && E_k^\prime \ar[rr] && F \ar[r] & 0 }.
\]
Since $Q$ and $F$ are torsion free sheaves, $E_k^\prime$ is torsion free.
Hence, $E_R^\prime$ is the extension of $E_K$ to $\mathcal X_R$,
corresponding to the vertex $[(E_R^\prime)_\xi]$ of $\mathfrak Q$.
On the other hand, we have the following exact sequence:
\[
\xymatrix@C=0.5cm{
  0 \ar[r] & \pi E_R \ar[rr] && E_R^\prime \ar[rr] && F\ar[r] & 0 }.
\]
Also, $E_R\overset{\pi}{\rightarrow}\pi E_R$ is an isomorphism. We get a homomorphism $E_R\overset{\pi}{\rightarrow}\pi E_R\hookrightarrow E^\prime_R$. Pulling back to $\mathcal X_k$, we get a homomorphism $E_k\rightarrow E_k^\prime$ and the image $Q^\prime$ of it is the subbundle corresponding to the edge $[E_\xi^\prime]-[E_\xi]$ at vertex $[E_\xi^\prime]$.
\end{proof}

The subbundle $F$ in Proposition \ref{prop beta bd}, is called the $\beta$-subbundle of the bundle $E$. For the convenience, in the following, the $\beta$-subbundle of $E$ should be denoted by $B$. Now assume that we are given a vertex $[E_\xi]$ of $\mathfrak Q$ such that the corresponding $E_k$ on
$\mathcal X_k$ is unstable. Let $B\subset E_k$ be the $\beta$-subbundle of $E_k$. Then, $\beta(B)>0$. By the Proposition \ref{prop b t c}, there is an edge in $\mathfrak Q$ at $[E_\xi]$
corresponding to $B$. Let $[E_\xi^{(1)}]$ be the vertex in $\mathfrak S$ determined by the edge, which corresponds to the subbundle $B$. Let $F_1\subseteq E_k^{(1)}$ be the image
of the canonical homomorphism $E_k\rightarrow E_k^{(1)}$(=the kernel of the homomorphism $E_k^{(1)}\rightarrow E_k$ ).

\begin{lemma}\label{lemm bb l}
If $G\subset E_k^{(1)}$ is a subbundle of $E_k^{(1)}$, then $\beta(G)\leq\beta(B)$, with equality
possible only  if $G\vee F_1=E_k^{(1)}$.
\end{lemma}

\begin{proof}
By the argument of Proposition \ref{prop b t c}, there are two exact sequences:
\[
\xymatrix@C=0.5cm{
  0 \ar[r] & B \ar[rr] && E_k \ar[rr] && F_1 \ar[r] & 0 };
\]
\[
\xymatrix@C=0.5cm{
  0 \ar[r] & F_1   \ar[rr] &&  E_k^{(1)}   \ar[rr] && B  \ar[r] & 0 }.
\]
If $G\subseteq F_1$, then there is a subbundle $W\subseteq E_k$,
such that
\[
\xymatrix@C=0.5cm{
  0 \ar[r] & B \ar[rr] && W  \ar[rr]  &&   G  \ar[r] & 0 }.
\]
Thus, $\beta(G)=\beta(W)-\beta(B)\leq 0$ (Proposition \ref{prop pro beta} and Proposition \ref{prop beta bd}). If $F_1\subset G$, then there is a subbundle $W^\prime\subseteq B$, such that
\[
\xymatrix@C=0.5cm{
  0 \ar[r] & F_1 \ar[rr]  && G  \ar[rr] &&  W^\prime \ar[r] & 0 }.
\]
So, $\beta(G)=\beta(F_1)+\beta(W^\prime)=\beta(W^\prime)-\beta(B)\leq 0$
($\beta(B)+\beta(F_1)=\beta(E_k)=0$ and Proposition \ref{prop pro beta}). For the other case, we have
$\beta(G)\leq\beta(G\vee F_1)+\beta(G\cap F_1)-\beta(F_1)\leq\beta(B)$ and equality
possible only  if $G\vee F_1=E_k^{(1)}$.
\end{proof}
Following \cite{lan}, we are now going to define a path $\mathcal P$ in $\mathfrak Q$
which starts at a given vertex $[E_\xi]$ such that the corresponding $E_k$ is unstable.
Let the succeeding vertex be the vertex determined by the edge corresponding to the $\beta$-subbundle $B$ of $E_k$.
If $\mathcal P$ reaches a vertex $[E_\xi^{(m)}]$ such that the corresponding bundle $E_k^{(m)}$ is semistable, then
the process stops automatically. If the path $\mathcal P$ never reaches a vertex corresponding
to a semistable reduction, then the process continuous indefinitely. In the following, We will show that the second alternative is impossible.\\
\hspace*{2em}Denote the $\beta$ subbundle of $E_k^{(m)}$ by $B^{(m)}$ and let $\beta_m=\beta(B^{(m)})$.
By Lemma \ref{lemm bb l}, $\beta_{m+1}\leq\beta_{m}$ and we must have $\beta_m>0$ unless $E_k^{(m)}$ is semistable. Thus, if the path $\mathcal P$ is continuous indefinitely we have $\beta_m=\beta_{m+1}=\cdots$ for sufficiently large $m$. Also, by Lemma \ref{lemm bb l}, for sufficiently large $m$, $B^{(m)}\vee F^{(m)}=E_k^{(m)}$, where $F^{(m)}=\text{Im}(E_k^{(m-1)}\rightarrow E_k^{(m)})$
($\text{Ker}(E_k^{(m)}\rightarrow E_k^{(m-1)})$). So, $\text{rank}(B^{(m)})+\text{rank}(F^{(m)})\geq r$. On the other hand,
$\text{rank}(B^{(m-1)})+\text{rank}(F^{(m)})=r$. Therefore, $\text{rank}(B^{(m)})\geq\text{rank}(B^{(m-1)})$, for sufficiently large $m$.
Since $\text{rank}(B^{(m)})\leq r$, we must have $\text{rank}(B^{(m)})=\text{rank}(B^{(m+1)})=\cdots$, for sufficiently large $m$.
Thus, $\text{rank}(B^{(m)})+\text{rank}(F^{(m)})= r$, $B^{(m)}\cap F^{(m)}=0$. Consequently, the canonical homomorphism
$E_k^{(m)}\rightarrow E_k^{(m-1)}$ induces an injection $B^{(m)}\hookrightarrow B^{(m-1)}$. Also, the canonical homomorphism
$E_k^{(m-1)}\rightarrow E_k^{(m)}$ induces an injection $F^{(m-1)}\hookrightarrow F^{(m)}$. Also, $\beta(B^{(m)})$ and $\text{rank}(B^{(m)})$ are both constant. It implies that $\beta(F^{(m)})=\beta(F^{(m+1)})=\cdots$, for $m$ sufficiently large.
\begin{lemma}
Let $R$ be a complete discrete valuation ring and $\mathcal P$ be an infinite path in $\mathfrak Q$, with vertices
$[E_\xi]$, $[E_\xi^{(1)}]$, $[E_\xi^{(2)}]$ $\cdots$. Let $F^{(m)}=\text{Im}(E_k^{(m+1)}\rightarrow E_k^{(m)})$. Assume that
$\text{rank}(F)=\text{rank}(F^{(1)})=\text{rank}(F^{(2)})=\text{rank}(F^{(3)})=\cdots=r$, the canonical homomorphism $E^{(m+1)}\rightarrow E^{(m)}$ induces injection $F^{(m+1)}\hookrightarrow F^{(m)}$, for each $m$, and $a_{n-1}(F)=a_{n-1}(F^{(1)})=a_{n-1}(F^{(2)})=\cdots$. Then $\beta(F)\leq 0$.
\end{lemma}
\begin{proof}
By the Lemma \ref{lemm ext 1}, there is a sequence of extensions of $E_K$ to $\mathcal X_R$, i.e
\[
\cdots\subset E^{(m)}\subset\cdots\subset E^{(1)}\subset E.
\]
Restricting the above inclusions to the special fiber $\mathcal X_k$, we get homomorphisms
\[
\cdots\rightarrow E^{(m)}_k\rightarrow\cdots\rightarrow E^{(1)}_k\rightarrow E_k
\]
and $F^{(m)}=\text{Im}(E_k^{(m+1)}\rightarrow E_k^{(m)})$, for each $m\geq 0$. Let $Q^{(m+1)}=\text{Ker}(E^{(m+1)}\rightarrow E^{(m)})$, for $m\geq 0$. Then the hypothesis that $F^{(m+1)}\hookrightarrow F^{(m)}$ is injective implies that $Q^{(m)}\cap F^{(m)}=(0)$. Let $\cdots\leftarrow E^{(m)}_k\leftarrow\cdots\leftarrow E^{(1)}_k\leftarrow E_k$ be the reverse homomorphisms.
By the Proposition \ref{prop b t c}, $Q^{(m)}=\text{Im}(E^{(m-1)}\rightarrow E^{(m)})$ and $F^{(m)}=\text{Ker}(E^{(m)}\rightarrow E^{(m+1)})$. Since $F^{(m)}\cap Q^{(m)}=(0)$, the induced map $Q^{(m)}\rightarrow Q^{(m+1)}$ is injective. By the exact sequence $\xymatrix@C=0.5cm{0 \ar[r] & F^{(m-1)} \ar[r] & E^{(m)} \ar[r] & Q^{(m)}\ar[r] & 0 }$, we have $a_{n-1}(F^{(m-1)})+a_{n-1}(Q^{(m)})=a_{n-1}(E^{(m)})=a_{n-1}(E_K)$, for $m\geq 1$. Since $a_{n-1}(F)=a_{n-1}(F^{(1)})=a_{n-1}(F^{(2)})=\cdots$, we have $a_{n-1}(Q^{(1)})=a_{n-1}(Q^{(2)})=a_{n-1}(Q^{(3)})=\cdots$.
Hence, the injections $Q^{(m)}\hookrightarrow Q^{(m+1)}$ are isomorphisms in codimension one. Also, $Q^{(m)**}$ are reflexive sheaves, then $Q^{(m)**}$ are determined by their restriction on the codimension one open substack. Thus, we have isomorphisms
\[
Q^{(1)**}\rightarrow Q^{(2)**}\rightarrow Q^{(3)**}\rightarrow Q^{(m)**}\rightarrow\cdots.
\]
So there is a sequence of inclusions:
\[
Q^{(1)}\hookrightarrow Q^{(2)}\hookrightarrow Q^{(3)}\hookrightarrow Q^{(m)}\hookrightarrow\cdots\hookrightarrow Q^{(1)**}.
\]
On the other hand, $Q^{(1)**}$ is a coherent sheaf on $\mathcal X_k$, it follows that
\[
Q^{(m)}\hookrightarrow Q^{(m+1)}\hookrightarrow Q^{(m+3)}\hookrightarrow\cdots.
\]
are isomorphisms, for sufficiently large $m$. Thus, we may assume without loss of generality that
\[
Q^{(1)}\hookrightarrow Q^{(2)}\hookrightarrow Q^{(3)}\hookrightarrow Q^{(m)}\hookrightarrow\cdots
\]
are isomorphisms. Also, we may assume that there is a subbundle $Q\subset E_k$ such that $Q\hookrightarrow Q^{(1)}$ is an isomorphism. Therefore, the exact sequence $\xymatrix@C=0.5cm{0 \ar[r] & F^{(m)} \ar[r] & E^{(m)}_k \ar[r]  & Q^{(m+1)}\ar[r] & 0 }$ splits, for each $m\geq 0$, i.e $E^{(m)}_k=F^{(m)}\oplus Q^{(m)}$. So, the exact sequence $\xymatrix@C=0.5cm{ 0 \ar[r] & Q^{(m+1)} \ar[r] & E_k^{(m+1)} \ar[r] &   F^{(m)} \ar[r] & 0 }$ yields $F^{(m+1)}\hookrightarrow F^{(m)}$ is an isomorphism, for each $m\geq 0$.\\
\hspace*{1em} Consider the completion $\hat{\mathcal X}_R$ of $\mathcal X_R$ with respect to the special fiber $\mathcal X_k$.
Let $\mathcal X_m=\mathcal X_R\times\text{Spec}(R/(\pi^m))$, for each $m\geq 0$. For a coherent sheaf $G$ on $\mathcal X_R$, we denote the restriction of $G$ to $\mathcal X_m$ by $G_m$. Following \cite{lan}, we will construct a coherent subsheaf $\hat{F}_R$ of $\hat{E}=\underset{{\longleftarrow}}{\lim}E_m$ on $\hat{\mathcal X}_R$. For each $m$, we will construct a coherent subsheaf $F_m$ of $E_m$ as following: Pulling back the inclusion $E^{(m)}\rightarrow E$ to $\mathcal X_m$, we get a homomorphism $E_m^{(m)}\rightarrow E_m$ and let $F_m$ be the image of this homomorphism. Let $j_{m,m^\prime}$ be the closed immersion $\mathcal X_{m^\prime}\hookrightarrow \mathcal X_m$, for $m^\prime\leq m$. Pulling back the homomorphism $E_m^{(m)}\twoheadrightarrow F_m\hookrightarrow E_m$ to $\mathcal X_{m^\prime}$, we get homomorphism $E_{m^\prime}^{(m)}\twoheadrightarrow j^*_{m,m^\prime}F_m\rightarrow E_{m^\prime}$, which fit into a commutative diagram:
\[
\xymatrix{
  E_{m^\prime}^{(m)} \ar[d]  \ar[r] & E_{m^\prime}^{(m^\prime)}\ar[d] \\
   j^*_{m,m^\prime}F_m \ar[r]  & E_{m^\prime}  }
\]
So, there is a natural homomorphism $j^*_{m,m^\prime}F_m\rightarrow F_{m^\prime}$. We can show this homomorphism is an isomorphism, step by step as the proof of Lemma 2 in \cite{lan}. Thus, we get an inverse system of sheaves $\{F_m\}$ and the inverse limit is a coherent subsheaf $\hat F_R$ of $\hat E$ on $\hat{\mathcal X}_R$.
By the Grothendieck's existence theorem for tame stacks in appendix A of \cite{av}, there exists a coherent subsheaf $F_R$ of $E$ such that $\hat{F}_R=\underset{{\longleftarrow}}{\lim}F_m$. Also, $j^*F_R=F$. Therefore, $a_{n-1}(F)=a_{n-1}(F_K)$, where $F_K=i^*F_R$. Since $E_K$ is semistable, we have $\beta(F)=\beta(F_K)\leq 0$.
\end{proof}

\begin{proof}[\textbf{The Proof of the second part in Theorem \ref{thm main 1}}]
For the case: $R$ is a complete discrete valuation ring, we have complete the proof of the Theorem \ref{thm main 1}.
As in \cite{lan}, the general case can be reduce to above case by considering the completion $\hat{R}$ of $R$. There is the following commutative diagram:
\[
\xymatrix{
  \mathcal X_k \ar[d]_{id} \ar[r]^{\hat{j}} & \mathcal X_{\hat{R}} \ar[d]_{p} & \mathcal X_{\hat{K}}\ar[l]_{\hat{i}}  \ar[d]^{p^\prime} \\
  \mathcal X_k \ar[r]^{j} & \mathcal X_R & \mathcal X_K \ar[l]_{i}  }
\]
Suppose that $E_K$ is a torsion free semistable sheaf on $\mathcal X_K$. Then the pullback ${p^\prime}^*E_K$ is a torsion free semistable sheaf on $\mathcal X_{\hat{K}}$(Proposition \ref{pro ex ss}), where $\hat{K}$ is the quotient field of $\hat{R}$. Denote the Bruhat-Tits complexes corresponding to $E_K$ and ${p^\prime}^*E_K$ by $\mathfrak Q_{1}$ and $\mathfrak Q_{2}$ respectively.
 For a vertex $[E_\xi]$ in the Bruhat-Tits complex $\mathfrak Q_1$, $[E_{\hat{R},\xi}]$ is the vertex in the Bruhat-Tits complex $\mathfrak Q_2$, where $E_{\hat{R},\xi}=E_\xi\underset{R}{\otimes}\hat{R}$. If $E$ is the torsion free sheaf on $\mathcal X_R$ corresponding to $[E_\xi]$, then $p^*E$ is the torsion free sheaf on $\mathcal X_{\hat{R}}$ corresponding to $[E_{\hat{R},\xi}]$. When $E_k$ is unstable, denote the $\beta$-subbundle of $E_k$ by $F_k$. By the Lemma \ref{prop b t c}, the edge $[E_\xi]-[E_{\xi}^\prime]$ in $\mathfrak Q_1$ corresponding to $F_k$ is constructed as following:
\[
\xymatrix@C=0.5cm{
  0 \ar[r] & F_k \ar[rr] && E_k \ar[rr] && Q_k \ar[r] & 0 }
\]
\[
\xymatrix@C=0.5cm{
  0 \ar[r] & E^\prime \ar[rr] && E \ar[rr] && Q_k \ar[r] & 0 }.
\]
And, the edge $[p^*E_\xi]-[p^*E_\xi^\prime]$ in $\mathfrak Q_2$ corresponding to $F_k$ is given by
\[
\xymatrix@C=0.5cm{
  0 \ar[r] & F_k \ar[rr] && E_k \ar[rr] && Q_k \ar[r] & 0 }
\]
\[
\xymatrix@C=0.5cm{
  0 \ar[r] & p^*E^\prime \ar[rr] && p^*E \ar[rr] && Q_k \ar[r] & 0 }.
\]
In $\mathfrak Q_2$, there is a finite path leading to a vertex whose corresponding torsion free sheaf on $\mathcal X_k$ is semistable, so it is in $\mathfrak Q_1$.
\end{proof}

\section*{Acknowledgement}
The author would like to thank Professor Jianxun Hu for his encouragement and help in the procedure of completing this paper. And also, thank Professor Yunfeng Jiang for carefully reading this paper and giving some valuable advices. We especially thank the referee for giving them very valuable advices which help them a lot to improve the presentations. This work was supported by the Fundamental Research Funds (34000-31610293) for the Central Universities, Sun Yat-sen University.



\begin{thebibliography}{99}

\bibitem {aov}D. Abramovich, M. Olsson and A. Vistoli, \textit{Tame stacks in positive characteristic}, Ann. Inst. Fourier(Grenoble), 58(4), 1057-1091 (2008).

\bibitem{av}D. Abramovich and A. Vistoli, \textit{Compactifying the Space of Stable Maps}, J. Am. Math. Soc., 15(1), 27-75 (2002).

\bibitem{sb}S. Brochard, \textit{Finiteness theorems for the Picard objects of an algebraic stack},
Adv. Math., 229, 1555-1585 (2012).

\bibitem{sb1}S. Brochard, \textit{Foncteur de Picard d'un champ alg\'ebrique}, Math. Ann. 343, 541-602 (2009).

\bibitem{bc}B. Conrad, \textit{Keel-Mori theorem via stacks}, unpublished manuscript.

\bibitem{dm}P. Deligne and D. Mumford, \textit{The irreducibility of the space of curves of given genus}, Publ. Math. IHES 36 (1969), 75-109.

\bibitem{Ha}R. Hartshorne, \textit{Algebraic Geometry}, GTM 52.(1977).


\bibitem{HJ}Y. H. Huang and Y. F. Jiang, \textit{The Geometry of Higgs bundles on Deligne-Mumford curves}, preprint.


\bibitem{hl}D. Huybrechts and M. Lehn, \textit{The Geometry of moduli spaces of sheaves}(Second Edition), Cambridge University Press.

\bibitem{ak}A. Kresch, \textit{On the geometry of Deligne-Mumford stacks}, In: Abramovich, D; Bertram, A; Katzarkov, L;
Pandharipande, R; Thaddeus, M. Algebraic Geometry: Seattle 2005. Providence, Rhode Island: Amer. Math. Soc. 259-271.

\bibitem{lan}S. G. Langton, \textit{Valuative Criteria for Families of Vector bundles on Algebraic varieties}, Annals of Mathematics Vol.101. No.1(1975). pp.88-110.

\bibitem{sl}S. Lang, \textit{Algebra (Revised Third Edition)}, GTM 211.(2002).


\bibitem{mm}M. Maruyama, \textit{Moduli of stable sheaves \uppercase\expandafter{\romannumeral 2}}, J. Math. Kyoto Univ. 18(1978), 557-614.

\bibitem{mr}V. Mehta and A. Ramanathan, \textit{An analogue of Langton's theorem on valuative criteria for vector bundles},
Proceedings of the Royal Society of Edinburgh. 96A, 39-45. 1984.

\bibitem{fn}F. Nironi, \textit{Moduli spaces of semistable sheaves on Projective Deligne-Mumford Stack},
 arXiv:0811.1949v2.

\bibitem{mo}M. Olsson, \textit{Algebraic spaces and stacks}, Amer. Math. Soc. Colloquium Publication, Volume 62, Amer. Math. Soc., Providence, RI, 2016.


\bibitem{os}M. Olsson and J. Starr, \textit{Quot functors for Deligne-Mumford stacks}, Comm. Algebra 31 (2003), no. 8, 4069-4096.

\bibitem{sp}C. Simpson, \textit{Moduli of representations of the fundamental group of a smooth projective variety\quad\uppercase\expandafter{\romannumeral 1}}, Publ. Math. IHES. 79(1994), 47-129.

\bibitem{bt}B. Toen, \textit{ Th\' eor\' emes de Riemann-Roch pour les champs de Deligne-Mumford}, K-Theory 18 (1999), 33-76.

\bibitem{vistoli}A. Vistoli, \textit{Intersection theory on algebraic stacks and on their moduli spaces}, Invent. math. 97 (1989), 613-670.

























\end{thebibliography}
\end{document}